\documentclass[12pt]{amsart}

\newtheorem{theorem}{Theorem}[section]

\newtheorem{corollary}[theorem]{Corollary}

\theoremstyle{definition}

\theoremstyle{remark}

\numberwithin{equation}{section}

\newcommand\lam{\lambda}

\newcommand\Del{\Delta}

\newcommand\Rc{\textup{Rc}}

\newcommand\ppt{\frac{\partial}{\partial t}}

\begin{document}

\title[Asymptotic behavior]{Some Asymptotic Behavior of the first Eigenvalue along the Ricci Flow}

\author[Jun Ling]{Jun \underline{LING}}
\address{Department of Mathematics, Utah Valley State College, Orem, Utah 84058}
\email{lingju@uvsc.edu}
\thanks{The author thank the Mathematical Sciences Research Institute at Berkeley for its
hospitality for program `The Geometric Evolution Equations and
Related Topics' and thank National Science Foundation for the
support offered.}



\date{October 8, 2007}



\begin{abstract}
We study some asymptotic behavior of the first nonzero eigenvalue
of the Lalacian along the normalized Ricci flow and give a direct
short proof for an asymptotic upper limit estimate.
\end{abstract}

\maketitle

\section{Introduction}\label{sec-intro}

The study of behavior of the eigenvalues of differential operators
along the flow of metrics is very active. We list a few such
studies as follows. Perelman \cite{perelman} proved the
monotonicity of the first eigenvalue of the operator $-\Del
+\frac14R$ along the Ricci flow by using his entropy and was then
able to rule out nontrivial steady or expanding breathers on
compact manifolds. X.~Cao \cite{caox} and  J.~F.~Li \cite{lijf}
studied the eigenvalues of $-\Del +\frac12R$ along the
unnormalized Ricci flow and gave some geometric applications. The
author \cite{lingj} studied the eigenvalues of Laplacians of the
normalized Ricci flow of metrics and gave a Faber-Krahn type of
comparison theorem and a sharp bound of the first nonzero
eigenvalue on compact 2-manifolds with negative Euler
Characteristic. In \cite{lingj2}, the author constructed a class
of monotonic quantities along the normalized Ricci flow. In this
short note, we study some asymptotic behavior of the first nonzero
eigenvalue of the Lalacian along the normalized Ricci flow and
give a direct short proof for an asymptotic upper limit estimate.

There are more developments to follow this.

\section{Results}\label{sec-results}
 We have the following results.
\begin{theorem}\label{thm1}
Let $M$ be an $n$-dimensional closed manifold. Let $g(t)$ be a
solution for $0\leq t<T$ ($T\leq \infty$) to the normalized Ricci
flow equation on $M$
\[
\ppt g=-2\Rc+\frac{2r}{n}g,
\]
where $\Rc$ is the Ricci tensor of the manifold $(M, g(t))$ and
$r$ the average of the scalar curvature $R$ of the manifold $(M,
g(t))$, and $g(t)$ starts with a Riemannian metric $g(0)$. Let
$\lam_{g(t)}$ be the first nonzero eigenvalue of the Laplacian
$\Del_{g(t)}$ of the metric $g(t)$. If $g(t)$ converges in all
$C^k(M)$-norm to a metric $g(T)$ as $t\rightarrow T$, then we
\[
\lam_{g(T)}\geq \limsup_{t\rightarrow T}\lam_{g(t)},
\]
where $\lam_{g(T)}$ is the first nonzero eigenvalue of the
Laplacian $\Del_{g(T)}$ of the limit metric $g(T)$.
\end{theorem}

\begin{proof}
Let $d\mu_t$ be the volume element of the metric $g(t)$, and
$d\mu_T$ the volume element of the metric $g(T)$.

For $\forall f\in C^{\infty}(M)$ with $\int_Mfd\mu_{T}=0$, let
\[
c_t=\int_M fd\mu_t.
\]
Then we have
\[
\lim_{t\rightarrow T}c_t=\int_M\lim_{t\rightarrow
T}fd\mu_t=\int_Mfd\mu_T=0.
\]
Since
\[
\begin{split}
\ppt(d\mu_t)
&=\frac{1}{2\sqrt{\det(g_{ij}(t))}}\left[\ppt\det(g_{ij}(t))\right]dx
\\
&=\frac{1}{2\det(g_{ij}(t))}\left[\ppt\det(g_{ij}(t))\right]d\mu_t
\\
&=\frac12g^{ji}(t)\left[\ppt g_{ij}(t)\right]d\mu_t
\\
&=(r-R)d\mu_t
\end{split}
\]
and
\[
r=\left\{\int_MRd\mu_t\right\}\Big/\left\{\int_Md\mu_t\right\},
\]
we have
\[
\ppt \int_Md\mu_t=\int_M(r-R)d\mu_t=0.
\]
Therefore $V=:\int_M d\mu_t=\int_M d\mu_T$ is constant independent
of $t$. Since the function
\[
h_t=:f-{c_t}/{V}
\]
satisfies the equation
 \[
 \int_Mh_td\mu_t=0,
 \]
the variational property of the first nonzero eigenvalue
(\cite{cha}, \cite{schoen-yau}) implies that
\[
\lam_{g(t)}\leq \frac{\int_M|\nabla^{[g(t)]}
h_t|^2d\mu_t}{\int_M(h_t)^2d\mu_t} =\frac{\int_M|\nabla^{[g(t)]}
f|^2d\mu_t}{\int_Mf^2d\mu_t-{(c_t)^2}/{V}},
\]
where $\nabla^{[g(t)]}$ is the gradient operator of the metric
$g(t)$.

Taking the upper limit as $t\rightarrow T$, we get
\[
\limsup_{t\rightarrow\infty}\lam_{g(t)}\leq
\frac{\int_M|\nabla^{[g(T)]} f|^2d\mu_T}{\int_Mf^2d\mu_T}.
 \]
By the variational property of the first nonzero eigenvalue again,
we have
\[
\lam_{g(T)}=\inf_{\stackrel{f\in C^{\infty}(M),}{\int_Mfd\mu_T=0}
}\frac{\int_M|\nabla^{[g(T)]}
f|^2d\mu_T}{\int_Mf^2d\mu_T}\geq\limsup_{t\rightarrow T}\lam(t).
\]
\end{proof}

\begin{corollary}\label{coro1}
Let $(M, g_0)$ be a closed Riemannian surface, or a closed
three-dimensional Riemannian manifold with positive Ricci
curvature, or a four-dimensional manifold with positive curvature
operator, and etc,  the normalized Ricci flow on $M$ with
$g(0)=g_0$, $g(\infty)$ the limit metric of $g(t)$ as
$t\rightarrow\infty$. Then we have
\[
\lam_{g(\infty)}\geq \limsup_{t\rightarrow \infty}\lam_{g(t)}.
\]
\end{corollary}
\begin{proof}
The conditions of Theorem \ref{thm1} are satisfied in each case,
by \cite{hamilton3}, \cite{hamilton}, and \cite{hamilton2},
respectively.
\end{proof}

\end{document}